\documentclass[a4paper,12pt,leqno]{article}
\usepackage{color}
\usepackage{amsmath, amssymb, amsthm,latexsym}
\usepackage{graphicx}
\usepackage[english]{babel}
\usepackage[pdftex,pagebackref=true,colorlinks=true,linkcolor=blue,citecolor=blue,urlcolor=blue]{hyperref}
\usepackage{currfile} %% To insert TeX file name.
\usepackage{fullpage}

\usepackage{version}
\excludeversion{add}

\theoremstyle{plain}
\newtheorem{theorem}{Theorem}[section]
\newtheorem{proposition}[theorem]{Proposition}
\newtheorem{lemma}[theorem]{Lemma}
\newtheorem{corollary}[theorem]{Corollary}

\newtheorem{state}[theorem]{Statement}

\theoremstyle{definition} %%\theoremstyle{remark}
\newtheorem{remark}[theorem]{Remark}
\newtheorem{remarks}[theorem]{Remarks}

\newtheorem{example}[theorem]{Example}

\theoremstyle{definition}
\newtheorem{definition}[theorem]{Definition}

\newtheorem{notation}[theorem]{Notation}
\numberwithin{figure}{section}

\newcommand{\R}{\mathbb{R}}

\newcommand{\cA}{\mathcal{A}}

\newcommand{\cF}{\mathcal{F}}

\newcommand{\noi}{\noindent}

\newcommand{\pf}{\noindent \textbf{Proof}.~}

\newcommand{\mf}{\mathfrak}
\newcommand{\mI}{\mathrm{I}}
\newcommand{\mJ}{\mathrm{J}}
\newcommand{\bc}{\bar{c}}
\newcommand{\hP}{\widehat{P}}

\newcommand{\vecp}[1]{\overset{\rightharpoonup}{#1}}

\newcommand{\sm}{\!\setminus\!}
\title{Sturm's theorem on the zeros of sums of eigenfunctions: Gelfand's strategy implemented}

\author{Pierre B\'erard and Bernard Helffer}

\date{\today}

\begin{document}
\maketitle

%\begin{abstract}
%\end{abstract}\medskip

\noindent Keywords: {Zeros of eigenfunction, Nodal domain, Courant nodal domain theorem,\\ Sturm theorem.}\\[5pt]
\noindent MSC~2010: {35P99, 35Q99, 58J50.}\\[5pt]
\noindent{\small Version:  {\currfilename}.}\\
%% Just for ``housekeeping purposes''.

%______________________________________________________

%% Local macros

\begin{abstract}In the second section ``Courant-Gelfand theorem'' of his last published paper (Topological properties of eigenoscillations in mathematical physics, Proc. Steklov Institute Math. 273 (2011) 25--34), Arnold recounts Gelfand's  strategy  to prove that the zeros of any linear combination of the $n$ first eigenfunctions of the Sturm-Liouville problem $$-\, y''(s) + q(x)\, y(x) = \lambda\, y(x) \mbox{ in } ]0,1[\,, \mbox{ with } y(0)=y(1)=0\,,$$
divide the interval into at most $n$ connected components, and concludes that ``the lack of a published formal text with a rigorous proof \dots is still distressing.''\\
Inspired by Quantum mechanics, Gelfand's strategy consists in replacing the ana\-lysis of linear combinations of the $n$ first eigenfunctions by that of their Slater determinant which is the first eigenfunction of the associated $n$-particle operator acting on Fermions.\\
In the present paper, we implement Gelfand's strategy, and give a complete proof of the above assertion. As a matter of fact,  refining Gelfand's strategy, we prove a stronger property taking the multiplicity of zeros into account, a result which actually goes back to Sturm (1836).
\end{abstract}

\section{Introduction}\label{S-intro}

On September 30, 1833, C.~Sturm\footnote{Jacques Charles Fran\c{c}ois \textsc{Sturm} (1803--1855)} presented a memoir on second order linear differential equations to the Paris Academy of Sciences. The main results are summarized in \cite{Sturm1833a,Sturm1833b}, and were later published in the first volume of Liouville's journal, \cite{Sturm1836a,Sturm1836b}. We refer to \cite{BH2017a} for more details. In this paper, we shall consider the following particular case.

\begin{theorem}[Sturm, 1836]\label{T-sturm}
Let $q$ be a smooth real valued function defined is a neighborhood of the interval $[0,1]$. The Dirichlet eigenvalue problem
\begin{equation}\label{E-intro-2}
\left\{
\begin{array}{l}
- \, y''(x) + q(x)\, y(x) = \lambda \, y(x) \text{~in~} ]0,1[\,,\\[5pt]
y(0)=y(1) = 0\,,
\end{array}%
\right.
\end{equation}
has the following properties.
\begin{enumerate}
  \item There exists an infinite sequence of (simple) eigenvalues
  $$\lambda_1 < \lambda_2 < \cdots \nearrow {}^{\infty},$$
  with an associated orthonormal family of eigenfunctions $\{h_j, j\ge 1\}$.
  \item For any $j \ge 1$, the eigenfunction $h_j$ has exactly $(j-1)$ zeros in the interval $]0,1[$.
  \item For any $1 \le m \le n$, let $U = \sum_{k=m}^{n} a_k\,h_k$ be any nontrivial real linear combination of eigenfunctions. Then,
      \begin{enumerate}
        \item $U$ has at most $(n-1)$ zeros in $]0,1[$, counted with multiplicities,
        \item $U$ changes sign at least $(m-1)$ times in $]0,1[$.
      \end{enumerate}
\end{enumerate}
\end{theorem}%

Sturm's motivations came from mathematical physics. He took a novel point of view, looking for qualitative behavior of solutions rather than for explicit solutions. To prove Assertions~1 and 2, he introduced the comparison and oscillation theorems which today bear his name.  Assertion~3  first appeared  as a corollary of Sturm's investigation of the evolution of zeros of a solution $u(t,x)$ of the associated heat equation, with initial condition $U$, as times goes to infinity (in direct line with his motivations).  Assertions~1 and 2 can be found in most textbooks on Sturm-Liouville theory. This is not the case for Assertion~3. In Section~\ref{S-pste}, we provide Liouville's proof,  which is based on the sole ordinary differential equation. We refer to \cite{BH2017a} for historical details.

\begin{remarks}\label{R-intro-2}
(i) In the framework of Fourier series, Assertion~3b is often referred to as the Sturm-Hurwitz theorem. See \cite{Ste2018}  for a quite recent qualitative version of this assertion.\\
 (ii) Sturm's theorem applies to more general operators, with more general boundary conditions; we refer to \cite{BH2017a} for more details.
\end{remarks}%

R.~Courant\footnote{Richard \textsc{Courant} (1888--1972).} partly generalized Assertion~2, in Sturm's theorem, to higher dimensions.

\begin{theorem}\label{T-cour}
Let $0 < \lambda_1 < \lambda_2 \le \lambda_3 \cdots \nearrow {}^{\infty}$ be the Dirichlet eigenvalues of $-\Delta$ in a bounded domain of $\R^d$, listed in nondecreasing order, with multiplicities. Let $u$ be any nontrivial eigenfunction associated with the eigenvalue $\lambda_n$, and let $\beta_0(u)$ denote the number of connected components of $ \Omega\sm u^{-1}(0)$ (\emph{nodal domains}). Then,
$$\beta_0(u) \le n\,.$$
\end{theorem}%

In a footnote of \cite[p.~454]{CH1953},  Courant and Hilbert make the following statement.

\begin{state}\label{T-gcour}
Any linear combination of the first $n$ eigenfunctions divides the domain, by means of its nodes, into no more than $n$ subdomains. See the G\"{o}ttingen dissertation of H.~Herrmann, Beitr\"{a}ge zur Theorie der Eigenwerten und Eigenfunktionen, 1932.
\end{state}%

In the literature, Statement~\ref{T-gcour} is referred to as the  ``Courant-Herrmann theorem'', ``Courant-Herrmann conjecture'',  ``Herrmann's theorem'', or ``Courant generalized theorem''. In \cite{BH2017b,BH2017c}, we call it the \emph{Extended Courant property.}

\begin{remarks}\label{R-intro-4}  Some remarks are in order. \vspace{-3mm}
\begin{enumerate}
  \item It is easy to see that Courant's upper bound is not sharp. This is indeed the case whenever the eigenvalue $\lambda_n$ is not simple. More generally, it can be shown that the number $\beta_0(u)$ is asymptotically smaller than $\gamma(n) \, n$ when $n$ tends to infinity, where $\gamma(n) < 1$ is a constant which only depends on the dimension $n$. It is interesting to investigate the eigenvalues for which Courant's upper bound is sharp, see the review article \cite{BoHe2017}. For this research topic, we also refer to the  surprising results in the recent paper \cite{JuZe2018}. \vspace{-3mm}
  \item In dimension greater than or equal to $2$, there is \emph{no general lower bound} for $\beta_0(u)$, except the trivial ones ($1$ for $\lambda_1$, and $2$ for $\lambda_k,k\ge 2$). Examples were first given by A.~Stern in her 1924 G\"{o}ttingen thesis, see \cite{BeHe2016}.
\end{enumerate}
\end{remarks}%

In the early 1970's, V.~Arnold\footnote{Vladimir Igorevich \textsc{Arnold} (1937-2010).} noticed that Statement~\ref{T-gcour},  would provide a partial answer to one of the problems formulated by D.~Hilbert\footnote{David \textsc{Hilbert} (1862--1943).}.

\begin{quotation}
\noi Citation from Arnold \cite[p.~27]{Arn2011}.\\[3pt]
I immediately deduced from the generalized Courant theorem [Statement~\ref{T-gcour}] new results in Hilbert's famous (16th) problem. \ldots\, And then it turned out that the results of the topology of algebraic curves that I had derived from the generalized Courant theorem contradict the results of quantum field theory. \ldots\, Hence, the statement of the generalized Courant theorem is not true (explicit counterexamples were soon produced by Viro). Courant died in 1972 and could not have known about this counterexample\footnote{As far as we know, the first paper of Arnold on this subject is \cite{Arn1973}, published in 1973.}.
\end{quotation}

Arnold was very much intrigued by Statement~\ref{T-gcour}, as is illustrated by \cite{Arn2011}, his last published paper, where he in particular relates a discussion with I.~Gelfand\footnote{Israel Moiseevich \textsc{Gelfand} (1913-2009).}, which we transcribe below, using Arnold's words, in the form of an imaginary dialog.\medskip

(Gelfand) \emph{I thought that, except for me, nobody paid attention to Cou\-rant's remarkable assertion. But I was so surprised that I delved into it and found a proof.}\smallskip

\noindent (Arnold is quite surprised, but does not have time to mention the counterexamples before Gelfand continues.)\smallskip

\emph{However, I could prove this theorem of Courant only for oscillations of one-dimensional media, where $m=1$.}\medskip

(Arnold) \emph{Where could I read it?}\medskip

(Gelfand) \emph{I never write proofs. I just discover new interesting things. Finding proofs (and writing articles) is up to my students.}\medskip

Arnold then recounts Gelfand's strategy to prove Statement~\ref{T-gcour} in the \emph{one-dimensional case}.

\begin{quotation}
\noi Quotations from \cite[Abstract and Section~2]{Arn2011}.\\[3pt]
\noi Nevertheless, the one-dimensional version of Courant's theorem is apparently valid. \ldots \, Gelfand's idea was to replace the analysis of the system of $n$ eigenfunctions of the one-particle quantum-mechanical problem by the analysis of the first eigenfunction of the $n$-particle problem (considering as particles, fermions rather than bosons). \ldots\medskip

\noi Unfortunately, [Gelfand's hints] do not yet provide a \emph{proof} for this generalized theorem: many facts are still to be proved. \ldots \medskip

\noi Gelfand did not publish anything concerning this: he only told me that he hoped his students would correct this drawback of his theory. \ldots \medskip

\noi Viktor Borisovich Lidskii told me that ``he knows how to prove all this''. \ldots \medskip

\noi Although [Lidskii's] arguments look convincing, the lack of a published formal text with a proof of the Courant-Gelfand theorem is still distressing.
\end{quotation}%

In \cite{Kuz2015}, Kuznetsov refers to Statement~\ref{T-gcour} as \emph{Herrmann's theorem}, and relates that Gelfand's approach \emph{so attracted Arnold that he included Herrmann's theorem for eigenfunctions of problem [\eqref{E-intro-2}]  together with Gelfand's hint into the 3rd Russian edition of his Ordinary Differential Equations}, see Problem~9 in the ``Supplementary problems'' at the end of \cite{Arn1992}.\medskip

More precisely, Arnold's Problem~9 proposes to prove the following statement, which is the one-dimensional analogue of Statement~\ref{T-gcour}.

\begin{state}\label{T-acour}
The zeros of any linear combination of the $n$ first eigenfunctions of the Sturm-Liouville problem \eqref{E-intro-2} divide the interval into at most $n$ connected components.
\end{state}%

This statement is equivalent to saying that any linear combination of the $n$ first eigenfunctions of \eqref{E-intro-2} has at most $(n-1)$ zeros in the open interval. This is a weak form of Sturm's upper bound, Assertion~3a in Theorem~\ref{T-sturm}. \medskip

In the present paper, we implement Gelfand's strategy to prove Statement~\ref{T-acour}, and we extend this strategy to take the multiplicities of zeros into account, and to prove Assertion~3a in Theorem~\ref{T-sturm}. Inspired by Quantum mechanics, Gelfand's strategy consists in replacing the analysis of linear combinations of the $n$ first eigenfunctions by that of their Slater determinant which is the first eigenfunction of the associated $n$-particle operator acting on Fermions. We give more details in Section~\ref{S-ho}. Note that Assertion~3b  can actually be deduced from Assertion~3a, see Section~\ref{S-pste}.\medskip

The paper is organized as follows. In Section~\ref{S-pste}, we give J.~Liouville's\footnote{Joseph \textsc{Liouville} (1809--1882).} 1836 proof of Assertion~3 in Theorem~\ref{T-sturm}. In Section~\ref{S-Nnot}, we introduce some notation. In Section~\ref{S-Nvdm}, we give preliminary results on Vandermonde polynomials, to be used later on. In Section~\ref{S-ho}, we explain Gelfand's strategy, and we apply it to a particular case, the harmonic oscillator.  Section~\ref{S-sl} is devoted to the proof of Assertion~3a in Theorem~\ref{T-sturm}, in the general case, following Gelfand's strategy: in Subsection~\ref{SS-sl2}, we prove Sturm's weak upper bound on the number of zeros of a linear combination of eigenfunctions, Statement~\ref{T-acour}, thus solving Problem~9 in \cite{Arn1992}; Sturm's strong upper bound is proved in Subsection~\ref{SS-sl4}.\bigskip

\textbf{Acknowledgements.} The authors would like to thank E.~Lieb and N.~Kuznetsov for useful comments on a first version of this paper.

\section{Liouville's proof of Sturm's theorem}\label{S-pste}

Assertions~1 and 2 in Theorem~\ref{T-sturm} are well-known, and can be found in many textbooks. This is not the case for Assertion~3. In this section, we give a short proof, based on the arguments of Liouville \cite{Liou1836b}, and Rayleigh\footnote{John William \textsc{Strutt},  Lord \textsc{Rayleigh} (1842--1919).} \cite[\S\, 142]{Ray1877}.\medskip

\emph{Proof of Assertion~3a.} Write equation~\eqref{E-intro-2} for $h_1$ and for $h_k$, multiply the first one by $h_k$, the second by $(-h_1)$ and add to obtain the relation
$$
\big( h_1 \, h_k' - h_1' \, h_k \big)' = (\lambda_1-\lambda_k)\, h_1 \, h_k\,.
$$
Multiply by $a_k$, and sum from $k=m$ to $k=n$ to obtain
\begin{equation}\label{E-pste-2}
\big( h_1 \, U' - h_1' \, U \big)' = h_1 \, U_1\,,
\end{equation}
where $U_1 = \sum_{k=m}^n (\lambda_1-\lambda_k) \, a_k \, h_k$.\\
 Integrating this relation from $0$ to $x$, and using the Dirichlet boundary condition, gives
$$
h_1(x) \, U'(x) - h_1'(x) \, U(x) = \int_{0}^x h_1(t) \, U_1(t) \, dt \,.
$$

Note that the left hand side can be rewritten as $h_1^2(x)\, \frac{d}{dx}\frac{U}{h_1}(x)$ in $]0,1[$. Count zeros with multiplicities. Assume that $U$ has $N$ zeros in $]0,1[$. Then so does $\frac{U}{h_1}$, so that, by Rolle's theorem, $\frac{d}{dx}\frac{U}{h_1}$ has \emph{a least} $(N-1)$ zeros in $]0,1[$. It follows that the function $x \mapsto \int_{0}^x h_1(t) \, U_1(t) \, dt$ has at least $(N-1)$ zeros in $]0,1[$. Note that it also vanishes at both $0$ and $1$ because the $h_j$ form an orthonormal family. By Rolle's theorem again, we conclude that its derivative, $h_1 \, U_1$, has at least $N$ zeros in $]0,1[$. Because $U$ and $U_1$ have the same form, we can repeat the argument, and conclude that, for any $\ell \ge 1$, the function $U_{\ell} = \sum_{k=m}^n (\lambda_1-\lambda_k)^{\ell} \, a_k \, h_k$ has \emph{at least} $N$ zeros in $]0,1[$. Letting $\ell$ tend to infinity, using the fact that the eigenvalues $\lambda_k$ are simple, and the fact that $h_n$ has $(n-1)$ zeros in $]0,1[$, it follows that $N \le (n-1)$. \medskip

\emph{Proof of Assertion~3b.} Assume that $U$ changes sign exactly $M$ times at the points $z_1 < \dots < z_M$ in the interval $]0,1[$, and that $M < (m-1)$, i.e., $M \le (m-2)$. Consider the function,
$$
V(x) :=
\begin{vmatrix}
h_1(z_1) &\dots &h_1(z_M) &h_1(x)\\
\vdots & &\vdots &\vdots \\
h_n(z_1) &\dots &h_n(z_M) &h_n(x)\\
\end{vmatrix}%
$$

It is easy to prove that the function $V$ is not identically zero (see Lemma~\ref{L-sl2-4}). It clearly vanishes at the points $z_j, 1 \le j \le M$, and it is a linear combination of the eigenfunctions $h_1,\dots,h_M$ (develop the determinant with respect to the last column). According to Assertion~3a in Theorem~\ref{T-sturm}, $V$ does not have any other zero, and each $z_j$ has order $1$, so that $V$ changes sign exactly at the points $z_j$. Since $M \le (m-2)$, the functions $U$ and $V$ are orthogonal, and their product $U\,V$ does not change sign in $]0,1[$. It follows that $U\,V$ vanishes identically, a contradiction. \hfill \qed

\begin{remark}\label{R-pste-2}
With the above notation, we can rewrite \eqref{E-pste-2} as
\begin{equation}\label{E-pste-4}
h_1 \, U_1 = h_1 \, U'' + (\lambda_1 - q) \, h_1\, U \,.
\end{equation}
A similar relation holds between $U_{\ell + 1}$ and $U_{\ell}$. Using these relations, and letting $\ell$ tend to infinity as in the preceding proof, we obtain the following lemma which is interesting in itself.
\end{remark}%

\begin{lemma}\label{L-pste-2}
The nonzero linear combination $U$ cannot vanish at infinite order at any point in $[0,1]$. In particular, its zeros are isolated.
\end{lemma}%

\section{Notation}\label{S-Nnot}

Let $n$ be an integer, $n\ge 1$, and $\mJ \subset \R$ an interval. Given $n$ points $x_1,\dots,x_n$ in $\mJ$, we denote the corresponding vector by $\vec{x} = (x_1,\dots,x_n) \in \mJ^n$. Generally speaking, we denote by $\vec{k} = (k_1, \cdots, k_n)$ a vector with positive integer entries. \medskip

We use the notation $\vecp{c} = (c_1, \dots, c_{n-1})$ for an $(n-1)$-vector with entries in $\mJ$.\medskip

Given $n$ real continuous functions $f_1, \dots, f_n$ defined on $\mJ$, we denote by $\vec{f}$ the vector-valued function $\big( f_1,\cdots,f_n\big)$, and we introduce the determinant
\begin{equation}\label{E-Nnot-2}
\left| \vec{f}(x_1) \dots \vec{f}(x_n)\right| :=
 \begin{vmatrix}
 f_1(x_1)&f_1(x_2)&\dots &f_1(x_n) \\
 f_2(x_1) &f_2(x_2)&\dots& f_2(x_n)\\
 \vdots&\vdots&\dots&\vdots \\
 f_n(x_1)&f_n(x_2)&\dots &f_n(x_n)
 \end{vmatrix}.
\end{equation}

Given a vector $\vec{b} = (b_1,\dots,b_n) \in \R^n$, we denote by
\begin{equation}\label{E-Nnot-4}
S_{\vec{b}}(x) = \sum_{j=1}^{n}b_j \, f_j(x)\,,
\end{equation}
the linear combination of $f_1, \dots, f_n$, with coefficients $b_j$'s.\medskip

Let $\vec{c} \in \mJ^n$ be a vector of the form
\begin{equation}\label{E-Nnot-10}
\vec{c} = (\bc_1,\dots,\bc_1,\bc_2,\dots,\bc_2,\,\dots\,, \bc_p,\dots,\bc_p ) \,,
\end{equation}
with $\bc_1$ repeated $k_1$ times, \ldots, $\bc_p$ repeated $k_p$ times, $1 \le p \le n$, $k_1 + \dots + k_p = n$, and with $\bc_1 < \bc_2 < \cdots < \bc_p$. \medskip

It will be convenient to relabel the variables $\vec{x} = (x_1,\dots,x_n)$ according to the structure of $\vec{c}$, as follows,
\begin{equation}\label{E-Nnot-12}
\vec{x}=(x_{1,1},\dots,x_{1,k_1},x_{2,1},\dots,x_{2,k_2},\dots,x_{p,1},
\dots,x_{p,k_p})\,,
\end{equation}
so that,
\begin{equation}\label{E-Nnot-12a}
\left\{
\begin{array}{l}
x_{1,1} = x_1, ~\dots, x_{1,k_1} = x_{k_1}\text{~and, for~}  2\le i \le p\,,\\[5pt]
x_{i,1} = x_{k_1+\cdots+k_{i-1}+1}, ~\dots, x_{i,k_i} = x_{k_1+\cdots+k_{i-1}+k_{i}}\,.
\end{array}
\right.
\end{equation}

In this case, we will also write the vector $\vec{x}$ as
\begin{equation}\label{E-Nnot-12b}
\vec{x} = \big(  x^{(1)},\dots,x^{(p)} \big)\,,
\end{equation}
with $x^{(i)} =(x_{i,1},\dots,x_{i,k_i})$, for $1 \le i \le p$.\medskip

We shall usually use both ways of labeling inside a formula, there should not be any confusion.\medskip

We introduce the real polynomials
\begin{equation}\label{E-Nnot-8q}
\left\{
\begin{array}{rl}
Q_1(x_1) &= 1\,, \text{~and, for~} n\ge 2,\\[5pt]
Q_n(x_1,\dots,x_n) &= \prod_{j=2}^{n}(x_1 - x_j)\,,
\end{array}%
\right.
\end{equation}
and
\begin{equation}\label{E-Nnot-8p}
\left\{
\begin{array}{rl}
P_1(x_1) &= 1\,, \text{~and, for~} n\ge 2,\\[5pt]
P_n(x_1,\dots,x_n) &= \prod_{1\le i < j \le n}(x_i - x_j) = \prod_{i=1}^{n} Q_{n+1-i}(x_i,\dots,x_n)\,.
\end{array}%
\right.
\end{equation}

\section{Vandermonde polynomials}\label{S-Nvdm}

\begin{lemma}\label{L-Nvdm-2}
The polynomial $P_n$, defined in \eqref{E-Nnot-8p}, is up to sign a Vandermonde\footnote{Alexandre Th\'{e}ophile \textsc{Vandermonde} (1735--1796).} determinant
\begin{equation}\label{E-Nvdm-2}
P_n(x_1,\dots,x_n) = (-1)^{\frac{n(n-1)}{2}}\, \begin{vmatrix}
1 & \dots & 1\\
x_1 & \dots & x_n\\
\vdots &  & \vdots\\
x_1^{n-1} & \dots & x_n^{n-1}
\end{vmatrix}.
\end{equation}
Furthermore,
\begin{enumerate}
  \item $P_n$ is anti-symmetric under the action of the group of permutations $\mf{s}_n$, and homogenous of degree $\frac{n(n-1)}{2}$.
  \item As a function of $x_1,\dots,x_n$, $P_n$ is harmonic, $\Delta P_n = 0$, and satisfies
  \begin{equation}\label{E-Nvdm-4}
  \partial^{n-1}_{x_n}\,\partial^{n-2}_{x_{n-1}}\cdots \, \partial^2_{x_3}\, \partial_{x_2}\,P_n = (n-1)!\,(n-2)!\,\dots\, 2!\,.
\end{equation}
\end{enumerate}
\end{lemma}

\proof The identity \eqref{E-Nvdm-2} is well-known, and readily implies Assertion~1. The polynomial $P_n$ being anti-symmetric, its Laplacian is also anti-symmetric, and hence, must be divi\-sible by $P_n$. Being of degree less than $P_n$, $\Delta P_n$ must be zero. The identity \eqref{E-Nvdm-4} follows immediately from the multi-linearity of the determinant, or by induction on $n$. \hfill \qed \medskip

\begin{notation}
When $\vec{x}=(x_1,\dots,x_n)$, we will also write $P_n(\vec{x})$ for $P_n(x_1,\dots,x_n)$ . We will denote by
$D_n(\partial_{\vec{x}})$ the differential operator which appears in \eqref{E-Nvdm-4}
\begin{equation}\label{E-Nvdm-4a}
D_n(\partial_{\vec{x}}) :=  \partial^{n-1}_{x_n}\,\partial^{n-2}_{x_{n-1}}\cdots \, \partial^2_{x_3}\, \partial_{x_2}\,,
\end{equation}
so that
\begin{equation}\label{E-Nvdm-4b}
D_n(\partial_{\vec{x}}) P_n(\vec{x}\,) = (n-1)!\,(n-2)!\,\dots\, 2!\,.
\end{equation}
\end{notation}%

\begin{notation}\label{N-omega}
In the sequel, we use $\omega(\vec{c},\vec{\xi}\,) $ as a generic notation for a function which depends on $\vec{c}, \vec{\xi}$, and tends to zero as $\vec{\xi}$ tends to zero.
\end{notation}%

\begin{lemma}\label{L-Nvdm-4}
Given $\vec{x} =(\vec{y},\vec{z}) \in \R^p\times \R^q$, the function
$$
\vec{x} \mapsto P_p(\vec{y})\, P_q(\vec{z})
$$
is harmonic as a function on $\R^{p+q}$.
\end{lemma}%

\medskip

We shall now describe the local behaviour of the harmonic polynomial $P_n$ near a point $\vec{c} \in \R^n$ at which it vanishes. We first treat two simple examples.

\begin{example}\label{Ex-Nvdm-1}
Let $n=5$, and $\vec{c}=(\bc_1,\bc_1,\bc_2,\bc_2,c_5)$, with $\bc_1 < \bc_2 < c_5$. Then, $P_5(\vec{c})=0$. Write $\vec{x} = \vec{c} + \vec{\xi}$. An easy computation gives,
\begin{equation}\label{E-Nvdm-Ex1}
P_5(\vec{c}+\vec{\xi})= P_2(\xi_1,\xi_2)\, P_2(\xi_3,\xi_4)\, \left\{ \rho(\vec{c}) + \omega(\vec{c},\vec{\xi}\,)\right\}\,,
\end{equation}
where $\rho(\vec{c}) = (\bc_1 - \bc_2)^4\,(\bc_1 - c_5)^2\,(\bc_2-c_5)^2$ is a nonzero constant.
\end{example}%

\begin{example}\label{Ex-Nvdm-2}
Let $n=5$. Let $\vec{c}=(\bc_1,\bc_1,\bc_1,c_4,c_5)$, with $\bc_1 < c_4 < c_5$. Then, $P_5(\vec{c})=0$. Write $\vec{x} = \vec{c} + \vec{\xi}$. An easy computation gives,
\begin{equation}\label{E-Nvdm-Ex2}
P_5(\vec{c}+\vec{\xi})= P_3(\xi_1,\xi_2,\xi_3)\, \left\{ \rho(\vec{c}) +  \omega(\vec{c},\vec{\xi}\,) \right\}\,,
\end{equation}
where $\rho(\vec{c}) = (\bc_1 - c_4)^3\,(\bc_1 - c_5)^3\,(c_4-c_5)$ is a nonzero constant.
\end{example}%

\begin{remark}\label{R-Ex12}
In both examples, the leading term on the right hand side of $P_n(\vec{c}+\vec{\xi}\,)$ is a homogeneous harmonic polynomial is some of the variables $\xi_j$'s, as we can expect from Bers's theorem, \cite{Ber1955}.  Furthermore, $\omega(\vec{c},\vec{\xi}\,)$ is actually a polynomial in the $(\xi_i-\xi_j)$'s, with coefficients depending on $\vec{c}$, and without constant term.
\end{remark}%

In the following lemma, we use both the standard coordinates names and their relabeling \eqref{E-Nnot-12}--\eqref{E-Nnot-12b}, for both variables $\vec{x}$ and $\vec{\xi}$.

\begin{lemma}\label{L-Nvdm-10}
Let $p$ be an integer, $1 \le p \le n$, and $(k_1,\dots,k_p)$ be a $p$-tuple of positive integers, such that $k_1 + \cdots + k_p =n$. Let $(\bc_1,\dots,\bc_p)$ be a $p$-tuple, such that $\bc_1 < \cdots < \bc_p$. Let $\vec{c}$ be the $n$-vector
\begin{equation}\label{E-Nvdm-10}
\vec{c} = ( \bc_1,\dots,\bc_1,\, \dots \,,\bc_p,\dots,\bc_p )\,,
\end{equation}
where each $\bc_j$ is repeated $k_j$ times, $1 \le j \le p$. Writing $\vec{x} = \vec{c} + \vec{\xi}$, and relabeling the coordinates of the vectors $\vec{x}$ and $\vec{\xi}$ as in \eqref{E-Nnot-12}--\eqref{E-Nnot-12b}, we have the following relation,
\begin{equation}\label{E-Nvdm-12}
P_n(\vec{c}+\vec{\xi}\,) = \rho(\vec{c}\,) \, P_{k_1}(\xi_{1,1},\dots,\xi_{1,k_1})\, \dots \, P_{k_p}(\xi_{p,1},\dots,\xi_{p,k_p})\left( 1 +  \omega (\vec{c},\vec{\xi}\,) \right),
\end{equation}
where $\rho(\vec{c}\,)$ is a nonzero constant depending only on $\vec{c}$, and where $ \omega (\vec{c},\vec{\xi}\,)$ is actually a polynomial in the variables $(\xi_i-\xi_j)$'s, with coefficients depending on the $c_j$'s,  without constant term.
\end{lemma}%

\proof  From the definition of $P_n$, and using the relabeling of the variables $\vec{x}$ and $\vec{\xi}$, as indicated in \eqref{E-Nnot-12}--\eqref{E-Nnot-12b}, we obtain the following relations.
\begin{equation}\label{E-Nvdm-16}
P_n(\vec{c}+\vec{\xi}\,) = \left( \prod_{i=1}^{k_1} Q_{n+1-i}(c_i+\xi_i,\dots,c_n+\xi_n)\right) \,
\prod_{i=k_1+1}^{n} Q_{n+1-i}(c_i+\xi_i,\dots,c_n+\xi_n)\,,
\end{equation}
\begin{equation}\label{E-Nvdm-18}
P_n(\vec{c}+\vec{\xi}\,) = \left( \prod_{i=1}^{k_1} Q_{n+1-i}(c_i+\xi_i,\dots,c_n+\xi_n)\right)\,
P_{n-k_1}(c_{2,1}+\xi_{2,1},\dots,c_{p,k_p}+\xi_{p,k_p})\,,
\end{equation}
Developing the factors $Q_{n+1-i}$ for $i\le k_1$, we obtain,
\begin{equation}\label{E-Nvdm-20}
\prod_{i=1}^{k_1} Q_{n+1-i}(c_i+\xi_i,\dots,c_n+\xi_n) = \rho_1(\vec{c}\,) \,
P_{k_1}(\xi_{1,1},\dots,\xi_{1,k_1})\,\left( 1 + \omega (\vec{c},\vec{\xi}\,) \right),
\end{equation}
where
\begin{equation}\label{E-Nvdm-22}
\rho_1(\vec{c}\,) = \left[ \left( \bc_1 - \bc_2 \right)^{k_2} \dots \left( \bc_1 - \bc_p \right)^{k_p}\right]^{k_1} \not = 0\,,
\end{equation}
and $\omega$ as in Notation~\ref{N-omega}.
%
%% \begin{equation}\label{E-Nvdm-22}
%% \rho_1(\vec{c}\,) = \left[ Q_{n+1-k_1}(\bc_1,c_{2,1},\dots,
%% c_{p,k_p}) \right]^{k_1} \not = 0\,,
%% \end{equation}
%
%and where
%
%\begin{equation}\label{E-Nvdm-24}
%\omega(\vec{c},\vec{\xi}\,)
%\end{equation}
%
%is a polynomial in the variables $(\xi_i-\xi_j)$'s, with coefficients depending on the %$c_j$'s,  without constant term.
%
%
Finally, we have
$$
P_n(\vec{c}+\vec{\xi}\,) = \rho_1(\vec{c}\,) \,
P_{k_1}(\xi_{1,1},\dots,\xi_{1,k_1})\, P_{n-k_1}(c_{2,1}+\xi_{2,1},\dots,c_{p,k_p}+\xi_{p,k_p})\,\left( 1 + \omega_1(\vec{c},\vec{\xi}\,) \right),
$$
or, more concisely,
\begin{equation}\label{E-Nvdm-25}
P_n(\vec{c}+\vec{\xi}\,) = \rho_1(\vec{c}\,) \,
P_{k_1}\big( \xi^{(1)} \big) \, P_{n-k_1}\big( c^{(2)} + \xi^{(2)},\dots,c^{(p)} +\xi^{(p)} \big) \,\left( 1 + \omega_1(\vec{c},\vec{\xi}\,)\right).
\end{equation}

\medskip

We can then apply the same kind of computation to the factor $P_{n-k_1}$, and repeat the operation until we finally obtain the desired formula, with
\begin{equation}\label{E-Nvdm-26}
\rho(\vec{c}\,) = \rho_1(\vec{c}\,) \, \cdots \, \rho_p(\vec{c}\,) \not = 0\,.
\end{equation}

\hfill \qed \medskip

We conclude this section with a technical lemma, which will play a key role later on.

\begin{lemma}[Division lemma]\label{L-DivL}
Let $P, Q$ be polynomials in $\R[X_1,\dots,X_n]$. Assume that $Q$ is harmonic and homogenous. If the set of real zeros of $Q$ is contained in the set of real zeros of $P$,
$$
\{x \in \R^n ~|~ Q(x)=0\} \subset \{x \in \R^n ~|~ P(x)=0\}\,,
$$
then $Q$ divides $P$, i.e. there exists $R$ in $\R[X_1,\dots,X_n]$ such that $P = Q R$.
\end{lemma}%

This lemma follows from Theorem~2 and Lemma~4 in \cite{Mur1964}. It is stated as Lemma~2.1 in \cite{LoMa2015}, with a proof given in \cite[Section~5.3]{LoMa2015}.

\section{Gelfand's strategy and the harmonic oscillator}\label{S-ho}

 In this section, we explain Gelfand's strategy to prove Statement~\ref{T-acour}, in the particular case of the harmonic oscillator. We also show how one can extend it to obtain a proof of Assertion~3a in Theorem~\ref{T-sturm}.\medskip

Let $\mf{H}^{(1)}$ denote the $1$-particle \emph{harmonic oscillator}
\begin{equation}\label{E-ho-2}
\mf{H}^{(1)}:=-\frac{d^2}{dx^2} + x^2
\end{equation}
on the line. The eigenvalues are given by $\{\lambda_n = 2n-1, n\ge 1\}$, they are simple, with associated orthonormal basis of eigenfunctions $\{h_n, n\ge 1\}$,
\begin{equation}\label{E-ho-4}
h_n(x) = \gamma_{n-1} \, H_{n-1}(x) \, \exp(-x^2/2)\,,
\end{equation}
where $H_m$ is the $m$-th Hermite polynomial, and $\gamma_m$ a normalizing constant \cite[Chap.~3]{Lau2014}. The polynomial $H_m(x)$ has degree $m$, with leading coefficient $2^m$, and satisfies the differential equation,
\begin{equation}\label{E-ho-6}
y''(x) - 2 x\, y'(x) + 2m\, y(x) = 0
\end{equation}
on the line $\R$. \medskip

We consider the $n$-particle Hamiltonian in $\R^n$,
\begin{equation}\label{E-ho-10}
\mf{H}^{(n)}:= \sum_{j=1}^n \left( - \frac{\partial^2}{\partial x_j^2} + x_j^2\right) = - \Delta + |\vec{x}|^2\,.
\end{equation}

\emph{Gelfand's strategy} is to look at $\mf{H}^{(n)}_F$, the operator $\mf{H}^{(n)}$ restricted to \emph{Fermions}, i.e., to functions which are anti-invariant under the action of the permutation group $\mf{s}_n$ on $\R^n$,
\begin{equation}\label{E-ho-12}
L_F^2(\R^n) = \left\lbrace f \in L^2(\R^n) ~|~ f\big( x_{\sigma(1)},\dots,x_{\sigma(n)}\big) = \varepsilon(\sigma) f(x_1,\dots,x_n), ~\forall \sigma \in \mf{s}_n \right\rbrace\,.
\end{equation}

Equivalently, we consider the Dirichlet realization $\mf{H}^{(n)}_F$ of $\mf{H}^{(n)}$ in
\begin{equation}\label{E-ho-14}
\Omega_n = \{(x_1,\dots,x_n) \in \R^n ~|~ x_1 < x_2 \cdots < x_n\}\,.
\end{equation}

Introduce the \emph{Slater}\footnote{John Clark \textsc{Slater} (1900--1976).}\emph{determinant}
\begin{equation}\label{E-ho-20}
\mf{S}_n(\vec{x}) = \begin{vmatrix}
                      h_1(x_1) & \dots & h_1(x_n) \\
                      \vdots &  & \vdots \\
                      h_n(x_1) & \dots & h_n(x_n) \\
                    \end{vmatrix} = A_n\,\exp(-|\vec{x}\,|^2/2)\, \begin{vmatrix}
                      H_0(x_1) & \dots & H_0(x_n) \\
                      \vdots &  & \vdots \\
                      H_{n-1}(x_1) & \dots & H_{n-1}(x_n) \\
                    \end{vmatrix}
                    .
\end{equation}

Using the properties of Hermite polynomials, we find that
\begin{equation}\label{E-ho-22}
\mf{S}_n(\vec{x}) = B_n \, \exp(-|\vec{x}\,|^2/2)\, P_n(\vec{x})\,.
\end{equation}
In the preceding equalities, $A_n$ and $B_n$ are nonzero constants depending only on $n$. \medskip

According to Arnold \cite[Section~2]{Arn2011}, Gelfand noticed the following two facts.

\begin{quotation}
\noi \textbf{A}.~The (antisymmetric) eigenfunction [$\mf{S}_n$] of the operator [$\mf{h}^{(n)}$] is the first eigenfunction for this operator (on functions satisfying the Dirichlet condition in the fundamental domain [$\Omega_n$]).\medskip

\noi \textbf{B}.~Choosing the locations [$(c_2,\dots,c_{n})$] of the other electrons (except for the first one), one can obtain any linear combination of the first $n$ eigenfunctions of the one-electron problem as a linear combination $[\mf{S}_n(x,c_2,\dots,c_n)]$ (up to multiplication by a nonzero constant).
\end{quotation}

Observe however that \textbf{B} is true only for linear combinations of the $n$ first eigenfunctions which have $(n-1)$ distinct zeros. \medskip

In the case of the harmonic oscillator, the proof of facts \textbf{A} and \textbf{B} is easy. More precisely, we have the following proposition which implies Statement~\ref{T-acour} in this particular case.\medskip

\begin{proposition}\label{P-ho-2}
Recall the notation $\vec{h}(c)=\big( h_1(c),\dots,h_n(c) \big)$.
\begin{enumerate}
  \item The function $\mf{S}_n(\vec{x}\,)$ is the first Dirichlet eigenfunction of $-\Delta + |\vec{x}\,|^2$ in $\Omega_n$.
  \item For any $\vecp{c}=(c_1,\dots,c_{n-1}) \in \Omega_{n-1}$, the vectors $\vec{h}(c_1),\dots,\vec{h}(c_{n-1})$, are linearly independent.
  \item Given $\vec{b} \in \R^n\sm\{0\}$, the linear combination
  $$
  S_{\vec{b}}(x) = \sum_{j=1}^n b_j h_j(x)
  $$
  has at most $(n-1)$ distinct zeros. Furthermore, if the function $S_{\vec{b}}$ has exactly $(n-1)$ distinct zeros $c_1 < c_2 < \dots < c_{n-1}$, then there exists a nonzero constant $C$ such that
$$
S_{\vec{b}}(x) = C \, \mf{S}_n(c_1,\dots,c_{n-1},x) \text{~for all~} x\in \R\,.
$$
\item The function $\mf{S}_n(c_1,\dots,c_{n-1},x)$ vanishes at order $1$ at each $c_j$, $1 \le j \le (n-1)$, and does not have any other zero.
\end{enumerate}
\end{proposition}%

\proof \emph{Assertion~1}. It is clear that $\mf{S}_n$ is an eigenfunction of $ - \Delta + |\vec{x}|^2$, and that it vanishes on $\partial \Omega_n$. From \eqref{E-Nvdm-2} and \eqref{E-ho-22}, we see that it does not vanish in $\Omega_n$, so that $\mf{S}_n$ must be the first Dirichlet eigenfunction for $-\Delta + |\vec{x}|^2$ in $\Omega_n$.\smallskip

\emph{Assertion~2}. If the vectors $\vec{h}(c_1),\dots,\vec{h}(c_{n-1})$, were dependent, $\mf{S}_n(c_1,\dots,c_{n-1},x)$ would be identically zero. Developing this determinant with respect to the last column, we would have $$\mf{S}_{n-1}(c_1,\dots,c_{n-1})\, h_n(x) + \dots \equiv 0\,.$$ This is impossible because the $h_j$'s are linearly independent and $\mf{S}_{n-1}(c_1,\dots,c_{n-1}) \not = 0\,$.\smallskip

\emph{Assertion~3}. Assume that $S_{\vec{b}}$ has at least $n$ distinct zeros $c_1 < \dots < c_n$. The $n$ components $b_j, 1 \le j \le n$ would satisfy a system of $n$ equations, whose determinant $\mf{S}_n(c_1,\dots,c_n)$ is nonzero. This would imply that $\vec{b}=\vec{0}$.  Assume that $S_{\vec{b}}$ has exactly $(n-1)$ zeros, $c_1 < \dots < c_{n-1}$. The function $x \mapsto \mf{S}_n(c_1,\dots,c_{n-1},x)$ can be written as a linear combination $S_{\vec{s}(\vecp{c})}(x)$, with coefficients $s_j(\vecp{c}), 1\le j\le n$ given by Slater like determinants. Both vectors $\vec{b}$ and $\vec{s}(\vecp{c})$ would then be orthogonal to the $(n-1)$ independent vectors $\vec{h}(c_1),\dots, \vec{h}(c_{n-1})$. This implies that there exists a nonzero constant $C$ such that $\vec{b} = C\, \vec{s}(\vecp{c})$.\smallskip

\emph{Assertion~4}. It suffices to consider the case of $c_1$.  Up to sign, we look at the local behavior of the function $x \mapsto \mf{S}(x,c_1,\dots,c_n)$  near $c_1$. Consider $\vec{c} =(c_1,c_1,c_2,\dots,c_{n-1})$, and write
$$
\mf{S}_n(\vec{c} + \vec{\xi}\,) = B_n \, \exp(- |\vec{c} + \vec{\xi}\,|^2/2) \,P_n(\vec{c} + \vec{\xi}\,) \,.
$$

Using Notation~\ref{N-omega} and Lemma~\ref{L-Nvdm-10}, we conclude that
$$
\mf{S}_n(\vec{c} + \vec{\xi}\,) = \alpha(\vec{c}\,)\, (\xi_1-\xi_2) \, \left( 1 + \omega(\vec{c},\vec{\xi}\,)\right) \,,
$$
for some nonzero constant $\alpha(\vec{c}\,)$ depending on $\vec{c}$.\medskip

It follows that
$$
\mf{S}_n(c_1+\xi,c_1,\dots,c_{n-1}) = \alpha(\vec{c}\,) \, \xi \, \left( 1 + \omega(\vec{c},\xi)\right) \,,
$$
so that this function vanishes precisely at order $1$ at $c_1$. \hfill \qed \medskip

\begin{remark}\label{R-ho-2}
It is standard in Quantum mechanics (except that the usual context for the one-particle Hamiltonian is a {3D}-space)  that the ground state energy of the $n$-particle Hamiltonian is the sum of the $n$ first eigenvalues of the one-particle Hamiltonian, a consequence of the Pauli\footnote{Wolfgang Ernst \textsc{Pauli} (1900--1958).} exclusion principle.  In a  context closer to our paper (see Section \ref{S-sl}), but with a different motivation, this sum associated with a one-particle Hamiltonian in an interval, and the properties of the corresponding  ground state,  are considered  in \cite{LiMa} at the beginning of the sixties.  Later on, this sum appears in  the celebrated Lieb-Thirring's inequality in connection with the analysis of the \emph{stability of matter} (see for example \cite{LiSe}) and references therein.
\end{remark}%

The following lemma allows us to extend  Gelfand's strategy in order to take care of the multiplicity of zeros,  and to achieve a proof of Sturm's upper bound.

\begin{lemma}\label{L-ho-2}
Let $\vecp{c} =(\bc_1,\dots,\bc_1,\,\dots\,,\bc_p,\dots,\bc_p)$, where $\bc_j$ is repeated $k_j$ times, with $\bc_1 < \dots < \bc_p$, and $k_1+\dots +k_p = n-1$. Let $\vec{k} =(k_1,\dots,k_p)$.  Define the function
\begin{equation}\label{E-ho-32}
\mf{S}_{\vec{k}}(x) = \left| \vec{h}(c_1)\dots\vec{h}^{(k_1-1)}(c_1)\,\dots\,
\vec{h}(c_p)\dots\vec{h}^{(k_p-1)}(c_p)\vec{h}(x) \right|,
\end{equation}
where $\vec{h}^{(m)}(x)$ is the vector $\big( h_1^{(m)}(x),\dots, h_n^{(m)}(x)\big)$, and where the superscript $(m)$ denotes the $m$-th derivative. \\ The function $\mf{S}_{\vec{k}}$ is not identically zero, and vanishes at exactly order $k_j$ at $\bc_j$. Furthermore, the vectors $\vec{h}(c_1),\dots,\vec{h}^{(k_1-1)}(c_1),\,\dots\,
,\vec{h}(c_p),\dots,\vec{h}^{(k_p-1)}(c_p)$, are linearly independent.
\end{lemma}%

\proof It suffices to consider the case of $\bc_1$. Clearly, $\mf{S}_{\vec{k}}$ vanishes at least at order $k_1$ at $\bc_1$.  It is sufficient to prove that the $k_1$-th derivative of this function does not vanish at $c_1$. We have
$$
\mf{S}_{\vec{k}}^{(k_1)}(x) = \pm\, \left| \vec{h}^{(k_1)}(x)\vec{h}(c_1)\dots\vec{h}^{(k_1-1)}(c_1)\vec{h}(c_2)\dots \vec{h}^{(k_2-1)}(c_2)\,\dots\,
\vec{h}(c_p)\dots\vec{h}^{(k_p-1)}(c_p) \right|.
$$

\emph{Claim: The value of this determinant at $x=\bc_1$ is different from zero}. Indeed, consider the vector $\vec{c}=(\bc_1,\dots,\bc_1,\,\dots\,,\bc_p,\dots,\bc_p)$, where $\bc_1$ is repeated $k_1+1$ times, and for $2\le j \le p$, $\bc_j$ is repeated $k_j$ times. Then $\mf{S}_{\vec{k}}^{(k_1)}(\bc_1)$ is a higher order derivative of $\mf{S}_n$ at $\vec{c}$. More precisely, using the relabeling of variables associated with $\vec{c}$, as given in \eqref{E-Nnot-12}--\eqref{E-Nnot-12b}, $\mf{S}_{\vec{k}}^{(k_1)}(\bc_1)$ is, up to sign, the derivative
$$
\left( \partial^{k_1}_{\xi_{1,k_1+1}}\dots \partial_{\xi_{1,2}}\right)\,
\left( \partial^{k_2-1}_{\xi_{2,k_2}}\dots \partial_{\xi_{2,2}}\right) \,\dots \,
\left( \partial^{k_p-1}_{\xi_{p,k_p}}\dots \partial_{\xi_{p,2}}\right) \, \mf{S}_n(\vec{c} + \vec{\xi}\,)\Big|_{\vec{\xi}=0}\,,
$$
or, using the notation \eqref{E-Nvdm-4a},
$$
D_{k_1}(\partial_{\xi^{(1)}})\dots D_{k_p}(\partial_{\xi^{(p)}})\, \mf{S}_n(\vec{c} + \vec{\xi}\,)\Big|_{\vec{\xi}=0}\,.
$$

The claim then follows from Lemma~\ref{L-Nvdm-2}, Equation~\eqref{E-Nvdm-4} and Lemma~\ref{L-Nvdm-10}, Equation~\eqref{E-Nvdm-12}.  The second assertion follows immediately. \hfill \qed \medskip

As a by product of the preceding proof, we have,

\begin{corollary}\label{C-ho-2}
Given, $p$, $1 \le p \le n$, let $k_1,\dots,k_p$ be $p$ positive integers such that $k_1+\cdots+k_p=n$. Let $\bc_1 < \dots < \bc_p$ be real numbers. Then, the determinant
\begin{equation}\label{E-ho-40}
\left| \vec{h}(c_1)\dots\vec{h}^{(k_1-1)}(c_1)\vec{h}(c_2)\dots \vec{h}^{(k_2-1)}(c_2)\,\dots\,\vec{h}(c_p)\dots\vec{h}^{(k_p-1)}(c_p) \right|
\end{equation}
is nonzero, so that the corresponding vectors are linearly independent.
\end{corollary}%

\begin{proposition}\label{P-ho-4}
For any $n\ge 1$, a nontrivial linear combination $S_{\vec{b}}$ of the eigenfunctions $h_1,\dots,h_n$ of the harmonic operator $\mf{H}^{(1)}$ has at most $(n-1)$ zeros on the real line, counted with multiplicities. Assume that $S_{\vec{b}}$ has $p$ zeros, $c_1 < \dots < c_p$ on the real line, with multiplicities $k_j$'s, such that $k_1+\cdots+k_p=n-1$. Then, there exists a nonzero constant $C$ such that
$$
S_{\vec{b}}(x) =C \, \left| \vec{h}(c_1)\dots\vec{h}^{(k_1-1)}(c_1)\,\dots\,
\vec{h}(c_p)\dots\vec{h}^{(k_p-1)}(c_p)\vec{h}(x) \right|.
$$
\end{proposition}%

\proof  The first assertion is Sturm's upper bound, Theorem~\ref{T-sturm}, in the particular case of the harmonic oscillator on the line.  The function $S_{\vec{b}}$  is a linear combination of the Hermite polynomials $H_0,\dots,H_{n-1}$, times the positive function $\exp(-|\vec{x}\,|^2/2)$. This immediately implies that the number of zeros of $S_{\vec{b}}$ on the real line, counted with multiplicities, is at most $(n-1)$.\medskip

 Here is a proof, \`{a} la Gelfand. \\
 Assume that a linear combination $S_{\vec{b}}$ has a least $n$ zeros on the real line, counted with multiplicities. From these zeros, one can determine some positive integer $p$, and sequences $\bc_1 < \dots < \bc_p$, $k_1,\dots,k_p$ satisfying the assumptions of  Corollary~\ref{C-ho-2}, and such that $S_{\vec{b}}$ vanishes at order (at least) $k_j$ at $\bc_j$, $ 1 \le j \le p$. This last condition implies that the $n$ entries of the vector $\vec{b}$ satisfy a system of $n$ equations, whose determinant is precisely $$|\vec{h}(c_1)\dots\vec{h}^{(k_1-1)}(c_1)\vec{h}(c_2)\dots \vec{h}^{(k_2-1)}(c_2)\,\dots\,\vec{h}(c_p)\dots\vec{h}^{(k_p-1)}(c_p)|.$$

Corollary~\ref{C-ho-2}  then implies that $\vec{b}=0$, so that a nontrivial linear combination $S_{\vec{b}}$ can have at most $(n-1)$ zeros on the real line, counted with multiplicities. \medskip

The second assertion is a consequence of (the proof of) Lemma~\ref{L-ho-2}.
\hfill \qed

\section{The Dirichlet Sturm-Liouville operator}\label{S-sl}

In this section, we show how Gelfand's strategy, see Section~\ref{S-ho}, can be applied to the general Dirichlet Sturm-Liouville problem~\eqref{E-intro-2}.

\subsection{Notation}\label{SS-sl1}

Let $q$ be a $C^{\infty}$ real function defined in a neighborhood of the interval $\mI := ]0,1[$. We  consider the $1$-particle operator
\begin{equation}\label{E-sl1-2}
\mf{h}^{(1)} := -\frac{d^2}{dx^2} + q(x)\,,
\end{equation}
and, more precisely, its Dirichlet realization in $\mI$, i.e. the Dirichlet boundary value problem
\begin{equation}\label{E-sl1-4}
\left\{
\begin{array}{l}
-\frac{d^2y}{dx^2} + q\,y = \lambda \, y\,,\\[5pt]
y(0) = y(1) = 0\,.
\end{array}%
\right.
\end{equation}

Let $\{(\lambda_j,h_j), j\ge 1\}$ be the eigenpairs of $\mf{h}^{(1)}$, with
\begin{equation}\label{E-sl1-6}
\lambda_1 < \lambda_2 < \lambda_3 < \cdots \,,
\end{equation}
and $\{h_j, j\ge 1\}$ an associated orthonormal basis of eigenfunctions.\medskip

We also consider the Dirichlet realization $\mf{h}^{(n)}$ of the $n$-particle operator in $\mI^n$,
\begin{equation}\label{E-sl1-12}
\mf{h}^{(n)} := - \sum_{j=1}^n \big( \frac{\partial^2}{\partial x_j^2} + q(x_j)\big) = - \Delta + Q\,,
\end{equation}
where $Q(x_1,\dots,x_n) = q(x_1) + \dots + q(x_n)$.\medskip

Denote by $\vec{k} = (k_1, \cdots, k_n)$ a vector with positive integer entries, and by $\vec{x} = (x_1,\cdots,x_n)$ a vector in $\mI^n$. The eigenpairs of $\mf{h}^{(n)}$ are the $(\Lambda_{\vec{k}},H_{\vec{k}})$, with
\begin{equation}\label{E-ls1-14}
\left\{
\begin{array}{l}
\Lambda_{\vec{k}} = \lambda_{k_1} + \cdots + \lambda_{k_n}\,, \text{~and~}\\[5pt]
H_{\vec{k}}(\vec{x}) = h_{k_1}(x_1) \cdots h_{k_n}(x_n)\,,
\end{array}%
\right.
\end{equation}
where $H_{\vec{k}}$ is seen as a function in $L^2(\mI^n,dx)$ identified with $ \widehat{\bigotimes} L^2(\mI,dx_j)$.\medskip

The symmetric group $\mf{s}_n$ acts on $\mI^n$ by $\sigma(\vec{x}) = (x_{\sigma(1)},\cdots,x_{\sigma(n)})$, if $\vec{x} = (x_1,\cdots,x_n)$. It consequently acts on $L^2(\mI^n)$, and on the functions $H_{\vec{k}}$ as well. A fundamental domain of the action of $\mf{s}_n$ on $\mI^n$ is the $n$-simplex
\begin{equation}\label{E-sl1-16}
\Omega_n^{\mI} := \{0 < x_1 < x_2 < \cdots < x_n < 1\}\,.
\end{equation}

In analogy with \eqref{E-ho-20}, we introduce the Slater determinant $\mf{S}_n$ defined by,
\begin{equation}\label{E-sl1-18}
 \mf{S}_n (x_1,\dots, x_n) =
 \begin{vmatrix}
 h_1(x_1)&h_1(x_2)&\dots &h_1(x_n) \\
 h_2(x_1) &h_2(x_2)&\dots& h_2(x_n)\\
 \vdots&\vdots&&\vdots \\
  h_n(x_1)&h_n(x_2)&\dots &h_n(x_n)
 \end{vmatrix}.
\end{equation}

Let $\vecp{c} =(c_1,\dots,c_{n-1}) \in \mI^{n-1}$. We consider the function $x \mapsto \mf{S}_n(c_1,\dots,c_{n-1},x)$. Developing the determinant with respect to the last column, we see that this function is a linear combination of the functions $h_1,\dots,h_n$, which we write as
\begin{equation}\label{E-sl1-22}
S_{s(\vecp{c})}(x) = \sum_{j=1}^{n}s_j(\vecp{c})\,h_j(x)
\end{equation}
where $s(\vecp{c}) = \left(s_1(\vecp{c}),\dots,s_n(\vecp{c}) \right)$, and
\begin{equation}\label{E-sl1-24}
s_j(\vecp{c})= s_j(c_1,\dots,c_{n-1}) = (-1)^{n+j} \,
\begin{vmatrix}
h_1(c_1) & \dots & h_1(c_{n-1})\\
\vdots &  & \vdots \\
h_{j-1}(c_1) & \dots & h_{j-1}(c_{n-1})\\
h_{j+1}(c_1) & \dots & h_{j+1}(c_{n-1})\\
\vdots &  & \vdots \\
h_n(c_1) & \dots & h_n(c_{n-1})\\
\end{vmatrix}
\end{equation}
so that $s(\vecp{c})$ is computed in terms of Slater determinants of size $(n-1)\times(n-1)$.

\subsection{Weak upper bound}\label{SS-sl2}

We now prove Statement~\ref{T-acour} using Gelfand's strategy, as explained in Section~\ref{S-ho}.

\begin{lemma}\label{L-sl2-4}
The function $\mf{S}_n$ is not identically zero.
\end{lemma}%

\proof The proof relies on the fact that the functions $h_j$, $1 \le j \le n$ are linearly independent. Clearly, $\mf{S}_1(x_1) = h_1(x_1) \not \equiv 0$. We now use induction on $n$. Assume that $\mf{S}_{n-1}(x_1,\dots,x_{n-1}) \not \equiv 0$. Develop the determinant $\mf{S}_n(x_1,\dots,x_n)$ with respect to the last column,
$$
\mf{S}_n(x_1,\dots,x_n) = \mf{S}_{n-1}(x_1,\dots,x_{n-1}) \, h_n(x) + \cdots \,.
$$
By the induction hypothesis, there exists $(x_1^0,\dots,x_{n-1}^0) \in \mI^{n-1}$, such that\\ $\mf{S}_{n-1}(x_1^0,\dots,x_{n-1}^0) \not = 0$. Then, $\mf{S}_n(x_1^0,\dots,x_{n-1}^0 ,x_n) \not \equiv 0$ because the $h_j$'s are linearly independent, and the lemma follows. \hfill \qed

\begin{lemma}\label{L-sl2-8}
The function $\mf{S}_n$ is the first Dirichlet eigenfunction of $\mf{h}^{(n)}$ in $\Omega_n^{\mI}$, with corresponding eigenvalue $\Lambda^{(n)}:=\lambda_1 + \cdots + \lambda_n$. In particular, the function $\mf{S}_n$ does not vanish in $\Omega_n^{\mI}$. More precisely, one can choose the signs of the functions $h_j$, $1 \le j \le n$, such that $\mf{S}_k$ is positive in $\Omega_k^{\mI}$ for $1 \le k \le n$. As a consequence, for any $c_1 < \dots < c_{n}$ in $\mI$, the vectors $\vec{h}(c_1),\dots,\vec{h}(c_n)$, are linearly independent.
\end{lemma}%

\pf An eigenfunction $\Psi$ of $\mf{h}_F^{(n)}$ is given by a (finite) linear combination $\Psi = \sum \alpha_{\vec{k}} H_{\vec{k}}$ of eigenfunctions of $\mf{h}^{(n)}$, such that the corresponding $\Lambda_{\vec{k}}$ are equal, and such that $\Psi$ is antisymmetric. If $\vec{k}=(k_1,\cdots, k_n)$ is such that $k_i = k_j$ for some pair $i\not = j$, using the permutation which exchanges $i$ and $j$, we see that the corresponding $\alpha_{\vec{k}}$ vanishes. It follows that the eigenvalues of $\mf{h}_F^{(n)}$ are the $\Lambda_{\vec{k}}$ such that the entries of $\vec{k}$ are all different. It then follows that the ground state energy of $\mf{h}_F^{(n)}$ is $\Lambda^{(n)}$.\medskip

It is clear that $\mathfrak S_n$ vanishes on $\partial \Omega^{\mI}_n$. Its restriction $\mathfrak S _{\Omega_n^{\mI}}$ to $\Omega_n^{\mI}$ satisfies the Dirichlet condition on $\partial \Omega_n^{\mI}$, and is an eigenfunction of $\mf{h}^{(n)}_F$ corresponding to $\Lambda^{(n)}$. Suppose that $ \mathfrak S _{\Omega_n^{\mI}}$ is not the ground state. Then, it has a nodal domain $\omega$ strictly included in $\Omega_n$. Define the function $U$ which is equal to $\mathfrak S _{\Omega_n^{\mI}}$ in $\omega$, and to $0$ elsewhere in $\mI^n$. It is clearly in $H_0^1(\Omega_n^{\mI})$. Using $\mathfrak s_n$, extend the function $U$ to a Fermi state $U_F$ on $\mI^n$. Its energy is $\Lambda^{(n)}$ which is the bottom of  the spectrum of $\mathfrak h^{(n)}_F$. It follows that $U_F$ is an eigenfunction of $\mathfrak h^{(n)}_F$, and a fortiori of $\mathfrak h^{(n)}$. This would imply that $\mf{S}_n$ is identically zero, a contradiction with Lemma~\ref{L-sl2-4}.\medskip

The fact that one can choose the $\mf{S}_n$ to be positive in $\Omega_n^{\mI}$ follows immediately. \medskip

If the vectors $\vec{h}(c_1),\dots,\vec{h}(c_n)$ were linearly dependent, the  function $\mf{S}_n$ would vanish at $(c_1,\dots,c_n) \in I$, a contradiction. \hfill \qed \medskip

%%% Simpler formulation.
%Since $\mf{S}_n$ is anti-invariant under the action of the permutation group %$\mf{s}_n$, $\mf{S}_n$ is an eigenfunction of $\mf{h}_F^{(n)}$, or
%equivalently of the Dirichlet realization of $- \Delta + Q$ in $\Omega_n^{\mI}$,
%with eigenvalue $\lambda_1+\dots+\lambda_n$. Since this eigenvalue is the
%least possible eigenvalue of $\mf{h}_F^{(n)}$ (a consequence of the %Pauli\footnote{Wolfgang Ernst \textsc{Pauli} (1900-1958).} exclusion
%principle), it follows that $\mf{S}_n$ is the first Dirichlet eigenfunction
%of $-\Delta + Q$ in $\Omega_n^{\mI}$, and hence does not vanish in
%this domain.

The following proposition provides a \emph{weak form} of Sturm's upper bound on the number of zeros of a linear combination of eigenfunctions of \eqref{E-sl1-4} (``weak'' in the sense that
 the multiplicities of zeros are not accounted for).

\begin{proposition}\label{P-sl2-swf}
Let $\vec{b} \in \R^n$, with $\vec{b} \not = \vec{0}$. Then, the linear combination $S_{\vec{b}}$ has a most $(n-1)$ distinct zeros in $\mI = ]0,1[$.
If $S_{\vec{b}}$ has exactly $(n-1)$ zeros in $\mI$, $c_1 < \cdots < c_{n-1}$, then there exists a nonzero constant $C$ such that
$$
S_{\vec{b}}(x) = C\, \mf{S}_n(c_1,\dots,c_{n-1},x)\,.
$$
Furthermore, each zero $c_j$ has order $1$.
\end{proposition}%

\proof Given $\vec{b}$, assume that $S_{\vec{b}}$ has at least $n$ distinct zeros $ c_1 < \dots < c_n$ in $\mI$. This means that the $n$ components $b_j, 1 \le j \le n$, satisfy the system of $n$ equations,
\begin{equation*}
\left\{
\begin{array}{l}
b_1 h_1(c_1) + \dots + b_n h_n(c_1) = 0,\\[5pt]
\cdots\\[5pt]
b_1 h_1(c_n) + \dots + b_n h_n(c_n) = 0.
\end{array}
\right.
\end{equation*}

By Lemma~\ref{L-sl2-8}, the determinant of this system is positive, and hence the unique possible solution is $\vec{0}$. This proves the first assertion. \medskip

Assume that $S_{\vec{b}}$ has precisely $(n-1)$ distinct zeros, $c_1 < \dots < c_{n-1}$, in $\mI$. By Lemma~\ref{L-sl2-8}, the vectors $\vec{h}(c_1),\dots,\vec{h}(c_{n-1})$, are linearly independent.  Then, $x \mapsto \mf{S}_n(c_1,\dots,c_{n-1},x)$ can be written as the linear combination $S_{\vec{s}(\vecp{c})}$, where the vector $\vec{s}(\vecp{c})$ is given by \eqref{E-sl1-24}. It follows that the vectors $\vec{b}$ and  $\vec{s}(\vecp{c})$ are both orthogonal to the family $\vec{h}(c_1),\dots,\vec{h}(c_{n-1})$, and must therefore be proportional. This proves the second assertion. \medskip

Assume that $x \mapsto \mf{S}_n(c_1,\dots,c_{n-1},x)$ vanishes at order at least $2$ at $c_1$. Then
$$
\frac{d}{dx}\Big|_{x=c_1}\mf{S}_n(x,c_1,c_2,\dots,c_{n-1}) =0.
$$
This implies that $ \frac{\partial \mf{S}_n}{\partial x_1}(c_1,c_1,c_2,\dots,c_{n-1}) = 0$, and hence that $\frac{\partial \mf{S}_n}{\partial \nu}(c_1,c_1,c_2, \dots,c_{n-1})$, where $\nu$ is the unit normal to the boundary $\partial \Omega_n^{\mI}$, which contradicts Hopf's lemma. This proves the last assertion, as well as the corollary.\hfill \qed \medskip

For completeness, we state the following immediate corollaries.

\begin{corollary}\label{C-sl2-swf1}
Given $c_1 < \cdots < c_{n-1}$ in $\mI$, the function
$$
x \mapsto \mf{S}_n(c_1,\dots,c_{n-1},x)\,,
$$
vanishes exactly at order $1$, changes sign at each $c_j$, and does not vanish elsewhere in $\mI$.
\end{corollary}%

\begin{corollary}\label{C-sl2-swf2}
Let $\vec{b} \in \R^n\sm\{0\}$. If the linear combination $S_{\vec{b}}$ has $k$ distinct zeros, and if one of the zeros has order at least $2$, then $k \le n-2$.
\end{corollary}%

\begin{remark}\label{R-sl2-swf}
Note that for $x \in ]c_j,c_{j+1}[$, $1 \le j \le n-1$,
$$
\mf{S}_n(c_1,\dots,c_{n-1},x) = (-1)^{n-1-j} \, \mf{S}_n(c_1,\dots,c_j,x,c_{j+1},\dots,c_{n-1}) \,,
$$
so that, according to Lemma~\ref{L-sl2-8}, it has the sign of $(-1)^{n-1-j}$. This also shows that this function of $x$ changes sign when $x$ passes one of the $c_j$'s.
\end{remark}%

%% Subsection on Gantmacher-Krein bound transferred to invisible Appendix
%% in ``version mode''

\subsection{Local behaviour of $\mf{S}_n$ near a zero}\label{SS-sl3}

We begin by treating two particular examples which are similar to Examples~\ref{Ex-Nvdm-1} and \ref{Ex-Nvdm-2}. We then deal with the general case. \medskip

Consider $\mf{S}_5$. Let $\vec{c} \in \partial \Omega_5^{\mI}$ be a boundary point. Write $\vec{x} = \vec{c} + \vec{\xi}$, with $\vec{\xi}$ close to $0$. The function $\mf{S}_5$ is an eigenfunction of the operator $-\Delta + Q$, and vanishes at the point $\vec{c} \in \mI^n$. By Bers' theorem \cite{Ber1955}, there exists a \emph{harmonic} homogeneous polynomial $\hP_k$, of degree $k$, such that
\begin{equation}\label{E-sl3-2}
\mf{S}_5(\vec{c}+\vec{\xi}\,) = \hP_k(\vec{\xi}\,) + \omega_{k+1}(\vec{c},\vec{\xi}\,)\,,
\end{equation}
where  $\omega_{k+1}(\vec{c},\vec{\xi}\,)$ is a function of $\vec{\xi}$, depending on $\vec{c}$, such that $\omega_{k+1}(\vec{c},t\vec{\xi}\,) = O(t^{k+1})$. Note that, for the time being, we have no a priori information on  the degree $k$.\medskip

\subsubsection{Example~1}\label{SSS-sl3-ex1} In this example, we take $\vec{c}=(\bc_1,\bc_1,\bc_2,\bc_2,c_5)$, with $\bc_1 < \bc_2 < c_5$. Call $\hP_k$ the polynomial given by \eqref{E-sl3-2} for this particular case.

\begin{lemma}\label{L-sl3-ex1}
The polynomial $\hP_k$  is given by
\begin{equation}\label{E-sl3-12}
\hP_k(\vec{\xi}\,) = \rho\,(\xi_1-\xi_2)(\xi_3-\xi_4)\,,
\end{equation}
where $\rho$ is a nonzero constant, and
\begin{equation}\label{E-sl3-12a}
\mf{S}_5(\vec{c}+\vec{\xi}\,) = \rho\, P_2(\xi_1,\xi_2)\, P_2(\xi_3,\xi_4)\, \big( 1 + \omega(\vec{c},\vec{\xi}\,) \big)\,,
\end{equation}
 where $\omega$ tends to zero when $\vec{\xi}$ tends to zero, see Notation~\ref{N-omega}.
\end{lemma}

\proof  According to \eqref{E-sl3-2}, we have
$$
\mf{S}_5(\bc_1+\xi_1,\bc_1+\xi_2, \bc_2+\xi_3,\bc_2+\xi_4,c_5+\xi_5) = \hP_k(\xi_1,\xi_2,\xi_3,\xi_4,\xi_5) + \omega_{k+1}(\vec{c},\vec{\xi}\,)\,.
$$
Using the anti-symmetry of $\mf{S}_5$, taking $\vec{\xi} = t\, \vec{\eta}$, using the fact that $\omega_{k+1}(\vec{c},t\, \vec{\eta}\,)$ is of order $k+1$, and letting $t$ tend to zero, we see that $\hP_k$ is anti-symmetric with respect to the pair $(\xi_1,\xi_2)$. A similar argument applies to the pair $(\xi_3,\xi_4)$.  This proves that
\begin{equation}\label{E-sl3-12ab}
\hP_k(\xi_1,\xi_2,\xi_3,\xi_4,\xi_5) = - \hP_k(\xi_2,\xi_1,\xi_3,\xi_4,\xi_5)
  = - \hP_k(\xi_1,\xi_2,\xi_4,\xi_3,\xi_5)\,,
\end{equation}
and hence, that $\hP_k(\xi_1,\xi_2,\xi_3,\xi_4,\xi_5)  = 0$ when $(\xi_1-\xi_2)(\xi_3-\xi_4) = 0$. \medskip

 We claim that the converse statement is true in a neighborhood of $0$. Indeed, assume that $\hP_k(\vec{\eta}\,)=0$, where $\eta_1 \not = \eta_2$ and $\eta_3 \not = \eta_4$.
  Using \eqref{E-sl3-12ab}, we can assume that $\eta_1 < \eta_2$ and $\eta_3 < \eta_4$. Because $\hP_k$ is a nonzero harmonic polynomial which vanishes at $\vec{\eta}$, in any neighborhood of $\vec{\eta}$, there exist points $\vec{\eta}^{\, \pm}$ such that $\hP_k(\vec{\eta}^{\,+}\,) \hP_k(\vec{\eta}^{\,-}\,) < 0$. For $t$ positive small enough, the function $\mf{S}_n(\vec{c}+t\, \vec{\eta}^{\, \pm}\,)$ has the sign of $\hP_k(\vec{c}+t\, \vec{\eta}^{\, \pm}\,)$, and this contradicts the fact that the function $\mf{S}_5$ is positive in $\Omega_n^{\mI}$.\medskip

 We have just proved that, in a neighborhood of zero, $\hP_k$ vanishes if and only if $(\xi_1-\xi_2)(\xi_3-\xi_4)$ vanishes.  The polynomials $\hP_k$ and $(\xi_1-\xi_2)(\xi_3-\xi_4)$ are both harmonic and homogeneous, and they have the same zero set in some neighborhood of zero. According to Lemma~\ref{L-DivL}, they divide each other, so that there exists a nonzero constant $\rho$ such that $\hP_k = \rho \, (\xi_1-\xi_2)(\xi_3-\xi_4)$. \hfill \qed

\subsubsection{Example~2}\label{SSS-sl3-ex2}  In this example, we choose $\vec{c}=(\bc_1,\bc_1,\bc_1,c_4,c_5)$, with $\bc_1 < c_4 < c_5$. Call $\hP_k$ the polynomial given by \eqref{E-sl3-2}.

\begin{lemma}\label{L-sl3-ex2}
The polynomial $\hP_k$ has the following properties. For any permutation
$\sigma \in \mf{s}_3(\xi_1,\xi_2,\xi_3)$, of the first three variables,
\begin{equation}\label{E-sl3-14}
\left\{
\begin{array}{l}
\hP_k(\xi_1,\xi_2,\xi_3,\xi_4,\xi_5) = \varepsilon(\sigma) \hP_k(\xi_{\sigma(1)},\xi_{\sigma(2)},\xi_{\sigma(3)},\xi_4,\xi_5)\,, \\[5pt]
 \hP_k = 0 \Leftrightarrow (\xi_1-\xi_2) (\xi_1-\xi_3)(\xi_2-\xi_3)=0 \,,\\[5pt]
\hP_k(\vec{\xi}\,) = \rho\, P_3(\xi_1,\xi_2,\xi_3)\,,
\end{array}
\right.
\end{equation}
where $\rho$ is a nonzero constant. This means that $\hP_k$ has degree $3$, and that
\begin{equation}\label{E-sl3-14a}
\mf{S}_5(\vec{c}+\vec{\xi}\,) = \rho\, P_3(\xi_1,\xi_2,\xi_3) \big( 1 + \omega(\vec{c},\vec{\xi}\,) \big)\,,
\end{equation}
where the function $\omega(\vec{c},\vec{\xi}\,)$ tends to zero when $\vec{\xi}$ tends to zero, see Notation~\ref{N-omega}.
\end{lemma}

\proof Similar to the previous proof. \hfill \qed

\subsubsection{General case}

Let $\vec{c} \in \partial \Omega_n^{\mI}$ be a boundary point, i.e. a point of the form $\vec{c} = (\bc_1,\dots,\bc_1,\, \dots\,,\bc_p, \dots,\bc_p)$, where $p$ is a positive integer, where $\bc_1 < \bc_2 < \dots <\bc_p$, are points in $\mI$, and where $\vec{c}$  is such that $\bc_j$ is repeated $k_j$ times, with $ k_1+\cdots +k_p =n$. \medskip

We write $\vec{x} = \vec{c} + \vec{\xi}$, with $\vec{\xi}$ close to $0$. The function $\mf{S}_n$ is an eigenfunction of the operator $-\Delta + Q$, and vanishes at the point $\vec{c} \in \mI^n$. By Bers's theorem \cite{Ber1955}, there exists a \emph{harmonic} homogeneous polynomial $\hP_k$, of degree $k$, such that
\begin{equation}\label{E-sl3-20}
\mf{S}_n(\vec{c}+\vec{\xi}\,) = \hP_k(\vec{\xi}\,) + \omega_{k+1}(\vec{c},\vec{\xi}\,)\,,
\end{equation}
where  the function $\omega_{k+1}(\vec{c},\vec{\xi}\,)$ is a function of $\vec{\xi}$, depending on $\vec{c}$, such that $\omega_{k+1}(\vec{c},t\vec{\xi}\,) = O(t^{k+1})$. Note that, for the time being, we have no a priori information on the degree $k$.\medskip

We relabel the coordinates of $\vec{\xi}$, according to \eqref{E-Nnot-12} -- \eqref{E-Nnot-12b}, and we write this vector as
\begin{equation}\label{E-sl3-22}
\vec{\xi} = \big( \xi^{(1)},\dots,\xi^{(p)}  \big)\,,
\end{equation}
where $\xi^{(j)} = (\xi_{j,1},\dots,\xi_{j,k_j})$.\medskip

The permutation group $\mf{s}_{k_j}$ acts by permuting the entries of $\xi^{(j)}$. Given $\sigma_j \in \mf{s}_{k_j}, 1\le j \le p$, we denote by $\sigma = (\sigma_1,\dots,\sigma_p) \in \mf{s}_{k_1} \times \dots \times \mf{s}_{k_p}$ the permutation in $\mf{s}_n$ which permutes the entries of $\xi^{(j)}$ by $\sigma_j$. \medskip

For the same vector $\vec{c}$, we look at the local behavior of  the Vandermonde polynomial  $P_n$, and we rewrite \eqref{E-Nvdm-12} as
\begin{equation}\label{E-sl3-26}
P_n(\vec{c}+\vec{\xi}\,) = \rho_1(\vec{c}\,) \, P_{k_1}\big( \xi^{(1)}\big) \cdots P_{k_p}\big(\xi^{(p)}\big)\left( 1 + \omega(\vec{c},\vec{\xi}\,) \right)\,.
\end{equation}

\begin{lemma}\label{L-sl3-gen}
The polynomial $\hP_k$ given by \eqref{E-sl3-20} has the following properties.
\begin{enumerate}
  \item  For any permutation $\sigma = (\sigma_1,\dots,\sigma_p) \in \mf{s}_{k_1} \times \dots \times \mf{s}_{k_p} \subset \mf{s}_n$,
  \begin{equation}\label{E-sl3-32a}
  \hP_k(\sigma\!\cdot\!\vec{\xi}\,) = \varepsilon(\sigma) \, \hP_k(\vec{\xi}\,)\,.
  \end{equation}
  \item The zero set of $\hP_k$ is characterized by
  \begin{equation}\label{E-sl3-32b}
  \hP_k(\vec{\xi}\,) = 0 \Leftrightarrow  \prod_{j=1}^p P_{k_j}\big( \xi^{(j)} \big) = 0 \,.
  \end{equation}
  \item There exists a nonzero constant $\rho(\vec{c}\,)$ such that
  \begin{equation}\label{E-sl3-32c}
  \hP_k(\vec{\xi}\,) = \rho(\vec{c}) \, P_{k_1}(\xi^{(1)})\dots P_{k_p}(\xi^{(p)})\,.
  \end{equation}
\end{enumerate}
This means that $\hP_k$ has degree $k= \sum_j \frac{k_j(k_j-1)}{2}$, and that
\begin{equation}\label{E-sl3-34}
\mf{S}_n(\vec{c}+\vec{\xi}\,) = \rho(\vec{c}\,)\, P_{k_1}(\xi^{(1)})\dots P_{k_p}(\xi^{(p)}) \big( 1 + \omega(\vec{c},\vec{\xi}\,) \big) \,,
\end{equation}
 where the function $\omega(\vec{c},\vec{\xi}\,)$ tends to zero when $\vec{\xi}$ tends to zero, see Notation~\ref{N-omega}.
\end{lemma}

\proof \emph{Assertion~1.} From the form of $\vec{c}$, and the definition of $\sigma = (\sigma_1,\dots,\sigma_p)$, we have the relations,
$$
\varepsilon(\sigma)\, \mf{S}_n(\vec{c} + t\vec{\xi}\,) = \mf{S}_n(\sigma\!\cdot\!(\vec{c}+t\vec{\xi}\,)) = \mf{S}_n(\vec{c}+t\sigma\!\cdot\!\vec{\xi}\,)\,.
$$
It follows that
$$
\hP_k(t\sigma\!\cdot\!\vec{\xi}\,) + \omega_{k+1}(\vec{c},t\sigma\!\cdot\!\vec{\xi}\,) = \varepsilon(\sigma) \, \big( \hP_k(t\vec{\xi}\,) + \omega_{k+1}(\vec{c},t\vec{\xi}\,)\big)\,.
$$
The assertion follows by dividing by $t$ and letting $t$ tend to zero.\medskip

\emph{Assertion~2.} The first assertion implies that the polynomial $\hP_{k}$  vanishes whenever the polynomial $\prod_{j=1}^p P_{k_j}\big( \xi^{(j)}\big)$ vanishes. Part $(\Leftarrow)$ of the second assertion follows.\smallskip

Assume that there exists some $\vec{\eta} =(\eta^{(1)},\dots,\eta^{(p)})$ such that \vspace{-2mm}
$$
\hP_k(\vec{\eta}\,) = 0 \text{~and~} \prod_{j=1}^p P_{k_j}\big( \eta^{(j)}\big) \not = 0.
$$
Since $\hP_{k}$ is harmonic, nonconstant, and vanishes at $\vec{\eta}$, it must change sign, and there exist $\vec{\eta}^{\,\pm}$ such that $\hP_k(\vec{\eta}^{\,+}) \hP_k(\vec{\eta}^{\,-}) < 0$. Using the first assertion and the properties of the Vandermonde polynomials, we see that one can choose $\vec{\eta}^{\,\pm} \in \Omega_n$, with $\Omega_n$ as in \eqref{E-ho-14}. It follows that for $t$ small enough, the vectors $\vec{c} + t \vec{\eta}^{\,\pm}$ are in $\Omega_n^{\mI}$, defined in \eqref{E-sl1-16}. For these vectors, one has
$$
\mf{S}_n(\vec{c} + t \vec{\eta}^{\,\pm}\,) = \hP_k(t \vec{\eta}^{\,\pm}) + \omega_{k+1}(\vec{c},t \vec{\eta}^{\,\pm})\,.
$$
This equality contradicts the fact that $\mf{S}_n$ is positive in $\Omega_n^{\mI}$. \medskip

\emph{Assertion~3.} Notice that the polynomials $\hP_k(\xi)$ and $\prod_{j=1}^p P_{k_j}\big( \xi^{(j)}\big)$  are both harmonic and homogeneous, with the same zero set in a neighborhood of $0$. We can then apply Lemma~\ref{L-DivL}, which implies that they divide each other, so that these polynomials must be proportional. The lemma is proved. \hfill \qed \medskip

As a consequence of the preceding lemma, we have,

\begin{corollary}\label{C-sl3-gen}
Let $\vec{c} \in \partial \Omega_n^{\mI}$ be as above. with the notation \eqref{E-Nvdm-4a}, we have the relations,
\begin{equation}\label{E-sl3-40}
D_{k_1}(\partial_{x^{(1)}})\cdots D_{k_p}(\partial_{x^{(p)}})\,\mf{S}_n(\vec{x}\,)\Big|_{\vec{x}=\vec{c}}\,
= D_{k_1}(\partial_{\xi^{(1)}})\cdots D_{k_p}(\partial_{\xi^{(p)}})\,\mf{S}_n(\vec{c}+\vec{\xi}\,)\Big|_{\vec{\xi} =0} \not = 0\,.
\end{equation}
\end{corollary}%

\subsection{Strong upper bound}\label{SS-sl4}

We can now prove Assertion~3a in Theorem~\ref{T-sturm}, using Gelfand's strategy, as explained in Section~\ref{S-ho}.

\begin{proposition}\label{P-sl4-2}
Let $\vec{b} \in \R^n\sm\{0\}$. Call $\bc_1 < \cdots < \bc_p$ the zeros of the linear combination $S_{\vec{b}}$ of the first $n$ eigenfunctions of problem \eqref{E-sl1-4}. Call $k_j$ the order of vanishing of $S_{\vec{b}}$ at $\bc_j$. Call $\vec{c}$ the vector $(\bc_1,\dots,\bc_1,\, \dots\,,\bc_p,\dots,\bc_p)$, where $c_j, 1 \le j \le p$ is repeated $k_j$ times. Then,
\begin{enumerate}
  \item $k_1 + \cdots + k_p \le (n-1)$,
  \item If $k_1 + \cdots + k_p = (n-1)$, then there exists a nonzero constant $C$ such that
  $$S_{\vec{b}} = C\, S_{\vec{s}(\vec{c})}\,,$$
where the linear combination $S_{\vec{s}(\vec{c})}$ is given by developing the determinant
\begin{equation}\label{E-sl4-2}
\left| \vec{h}(c_1)\dots\vec{h}^{(k_1-1)}(c_1)\, \dots \,   \vec{h}(c_p)\dots\vec{h}^{(k_p-1)}(c_p) \vec{h}(x) \right|,
\end{equation}
and where $\vec{h}^{(m)}(a)$ is the vector $\big(h_1^{(m)}(a),\dots,h_n^{(m)}(a)\big)$ of the $m$th derivatives of the $h_j$'s evaluated at the point $a$.
\end{enumerate}

\end{proposition}%

\proof \emph{Assertion~1.} Assume that $k_1+\cdots+k_p \ge n$.  This implies that the coefficients $b_1,\dots,b_n$, satisfy the system of $n$ equations,
$$
(b_1,\dots,b_n) \left( \vec{h}(c_1)\dots\vec{h}^{(k_1-1)}(c_1)\, \dots \,   \vec{h}(c_p)\dots\vec{h}^{(k_p-1)}(c_p) \right) = 0
$$
where the left hand side is the product of the row matrix $(b_1,\dots,b_n)$ by the $n\times n$ matrix
$$
\left( \vec{h}(c_1)\dots\vec{h}^{(k_1-1)}(c_1)\, \dots \,   \vec{h}(c_p)\dots\vec{h}^{(k_p-1)}(c_p) \right)\,.
$$
Using \eqref{E-sl3-40}, we see that the determinant of the latter matrix is nonzero. This implies that $\vec{b}=0$, a contradiction.\medskip

\emph{Assertion~2.} Using \eqref{E-sl3-40} again (with $n-1$ instead of $n$), we see that the coefficient of $h_n(x)$ in the linear combination $S_{\vec{s}(\vec{c})}$ is nonzero, so that $S_{\vec{s}(\vec{c})}$ is not identically zero. It follows that the family of $(n-1)$ vectors $\cF:= \left\lbrace \vec{h}(c_1),\dots,\vec{h}^{(k_1-1)}(c_1),\, \dots \,,   \vec{h}(c_p),\dots ,\vec{h}^{(k_p-1)}(c_p)\right\rbrace$ is free. Both functions $S_{\vec{b}}$ and $S_{\vec{s}(\vec{c})}$ vanish at order $k_j$ at $\bc_j$, for $1 \le j \le p$. This means that the vectors $\vec{b}$ and $\vec{s}(\vec{c})$ are both orthogonal to $\cF$, which implies that they are proportional. The proposition is proved. \hfill \qed \medskip

\begin{remark}\label{R-gk}
In this paper, we have considered a Dirichlet Sturm-Liouville problem with smooth coefficients. In less regular cases, one can still improve Statement~\ref{T-acour} by introducing the number $N_{\vec{b}}$ of nodes of $S_{\vec{b}}$ (zeros at which the function changes sign), and the number $A_{\vec{b}}$ of anti-nodes (zeros at which the function retains its sign). Then, $N_{\vec{b}} + 2 A_{\vec{b}} \le n-1$. This result is stated in \cite[p.~275]{GZ2003}, and proved in \cite[Chap.~{III}.5]{GaKr2002} in the more general framework of Chebyshev systems of continuous functions.
\end{remark}%

\begin{add}%% Version addendum
\newpage
\appendix

\section{The upper bound of Gantmacher and Krein}\label{A-gk}

In this Appendix, we describe an upper bound which is stronger than Statement~\ref{T-acour}, and weaker than Theorem~\ref{T-sturm}, Assertion~(3a). It is mentioned in \cite{GZ2003} in the particular case of the Dirichlet Sturm-Liouville eigenvalue problem, and proved in the book by \cite{GaKr2002} Gantmacher and Krein, in the more general framework of Chebyshev systems of continuous functions.

\begin{definition}\label{D-sl2-swf}
Let $S$ be a continuous function in $\mI$. Let $c \in \mI$ be a zero of $S$. Following \cite[Chap.~{III}.5]{GaKr2002}, call $c$ a \emph{node} of $S$, if for any $\varepsilon > 0$,  small enough, there exists some $x_{\varepsilon}^{\pm}$ such that $c-\varepsilon < x_{\varepsilon}^{-} < c < x_{\varepsilon}^{+} < c+\varepsilon$, with $S(x_{\varepsilon}^{-})\, S(x_{\varepsilon}^{+}) < 0$; call $c$ an \emph{antinode} of $S$,  if for any $\varepsilon > 0$, small enough, $S$ does not change sign in $]c-\varepsilon,c+\varepsilon[$, and does not vanish identically in $]c-\varepsilon,c[$ and in $]c,c+\varepsilon[$.
\end{definition}%

This definition applies to any continuous function $S$. An isolated zero of $S$ is either a node or and antinode. In the particular case of the Sturm-Liouville eigenvalue problem with smooth coefficients, according to Lemma~\ref{L-pste-2}, any zero of a nontrivial $S_{\vec{b}}$ is isolated. It is a node (resp.  an antinode), if and only if the first nonzero coefficient in  the Taylor series of $S$ at $c\,$ is odd (resp. even). The following proposition is a consequence of \cite[Chap.~{III}.5]{GaKr2002}. We  reproduce the proof for completeness.

\begin{proposition}\label{P-sl2-gm}
Let $\vec{b} \in \R^n\sm\{0\}$. Let $N_{\vec{b}}$ be the number of nodes of $S_{\vec{b}}$, resp. let $A_{\vec{b}}$ be the number of antinodes. Then,
$$
N_{\vec{b}} + 2\, A_{\vec{b}} \le n-1\,.
$$
\end{proposition}%

\proof In this proof, we use the notation $S$ for $S_{\vec{b}}$. We already know, Proposition~\ref{P-sl2-swf}, that $S$ has at most $(n-1)$ distinct zeros. \medskip

We say that a set $z_1 < z_2 < \dots < z_s$ has the property $\cA$ with respect to $S$ (alternating property) if there exists $\kappa \in \{0,1\}$ such that for any
$k \in \{1,\dots,s\}$, $(-1)^{k+\kappa}\, S(z_k) \ge 0$.\medskip

\begin{lemma}\label{L-sl2-gm1}
Assume that $\{z_1 < z_2 < \dots < z_s\}$ has the property $\cA$ with respect to $S$. If $\xi \not \in \{z_1 < z_2 < \dots < z_s\}$ is an antinode of $S$, then, for $\varepsilon$ small enough, one of the sets $\{z_1 < z_2 < \dots < z_s\} \cup \{\xi-\varepsilon,\xi\}$ or $\{z_1 < z_2 < \dots < z_s\} \cup \{\xi,\xi+\varepsilon\}$, properly reordered, has the property $\cA$, with $S(\xi - \varepsilon) \, S(\xi + \varepsilon) > 0$.
\end{lemma}%

\emph{Proof of the lemma}. We examine the case in which there exists $1 < j < s-1$, such that $z_j < \xi < z_{j+1}$. We choose $\varepsilon$ such that
$$
z_j < \xi-\varepsilon < \xi < \xi+\varepsilon < z_{j+1}\,,
$$
with $S(\xi \pm \varepsilon) \not = 0$.\medskip

We know that $(-1)^{\kappa+j}S(z_j) \ge 0$,  and that $S(\xi)=0$. We have two cases,
\begin{itemize}
\item if $(-1)^{\kappa+j+1} S(\xi-\varepsilon) > 0$, then
$(-1)^{\kappa+j+2} S(\xi) \ge 0$, and $(-1)^{\kappa+j+3} S(z_{j+1}) \ge 0$,
\item if $(-1)^{\kappa+j+1} S(\xi-\varepsilon) < 0$, then
$(-1)^{\kappa+j+1} S(\xi) \ge 0$, and $(-1)^{\kappa+j+2} S(\xi+\varepsilon) \ge 0$,
\end{itemize}

From the set $\{z_1,\dots, z_s,\xi-\varepsilon,\xi,\xi+\varepsilon\}$, we construct
an ordered list $\{z'_1 < \dots < z'_{s+2}\}$ as follows
\begin{itemize}
\item for $k \le j$, $z'_k = z_k$, and for $k \ge j+3$, $z'_k = z_{k+2}$,
and we choose $z'_{j+1} < z'_{j+2}$ in the interval $]z_j,z_{j+1}[$, as follows:
\item if $(-1)^{\kappa+j+1}S(\xi-\varepsilon) > 0$, we choose
$$
\left\{
\begin{array}{l}
z'_{j+1} = \xi-\varepsilon\,, \text{~so that~} (-1)^{\kappa+j+1} S(z'_{j+1}) > 0,\\
z'_{j+2} = \xi\,, \text{~so that~} (-1)^{\kappa+j+2} S(z'_{j+2}) = 0,\\
\end{array}
\right.
$$
\item if $(-1)^{\kappa+j+1}S(\xi-\varepsilon) < 0$, then
$(-1)^{\kappa+j+2}S(\xi+\varepsilon) > 0$, and we choose
$$
\left\{
\begin{array}{l}
z'_{j+1} = \xi \,, \text{~so that~} (-1)^{\kappa+j+1} S(z'_{j+1}) = 0\,,\\
z'_{j+2} = \xi+\varepsilon\,, \text{~so that~} (-1)^{\kappa+j+2} S(z'_{j+2}) > 0\,.\\
\end{array}
\right.
$$
\end{itemize}

The case $\xi < z_1$ or $\xi > z_s$ are dealt with similarly. This proves the lemma. \hfill \qed \medskip

\emph{Proof of the proposition continued.} Call $z_1 < \dots < z_p$, $p = N_{\vec{b}}$, the nodes of $S$. Then, one can choose numbers $\alpha_j$ such that,
$$
\alpha_1 < z_1 < \alpha_2 < z_2 < \cdots < \alpha_p < z_p < \alpha_{p+1}\,,
$$
and some $\kappa \in \{0,1\}$ such that $(-1)^{\kappa+j}S(\alpha_j) > 0\,$. \medskip

By applying Lemma~\ref{L-sl2-gm1} recursively for the $q = A_{\vec{b}}$ antinodes, we obtain a set
$$
\beta_1 < \beta_2 < \cdots < \beta_{p+2q+1}
$$
such that $(-1)^{\kappa+j}S(\beta_j) \ge 0$ (Property $\cA$). \medskip

Assume that $p+2q > n-1$, i.e. $p+2q+1 \ge n+1$. Consider the vector $\vec{H} =(h_1,\dots,h_n,S)$. Then, the determinant
$$
|\vec{H}(\beta_1)\dots\vec{H}(\beta_{n+1})| \text{~is identically~} 0,
$$
because $S$ is a linear combination of $h_1,\dots,h_n$. Developing this determinant with respect to the last row, we find that
$$
0 \equiv \sum_{k=1}^n (-1)^{n+1+k}S(\beta_k) \, |\vec{h}(\beta_1)\dots\vec{h}(\beta_k-1)\vec{h}(\beta_k+1)\dots\vec{h}(\beta_n)|\,.
$$

For each $k$, we have $(-1)^{n+1+k}S(\beta_k) \ge 0$ by construction of the set $\{\beta_j, 1\le j \le n+1\}$, and
$$
|\vec{h}(\beta_1)\dots\vec{h}(\beta_{k-1})\vec{h}(\beta_{k+1})\dots\vec{h}(\beta_n)|
 = \mf{S}_n(\beta_1\dots\beta_{k-1},\beta_{k+1},\dots,\beta_{n+1}) > 0
$$
according to Lemma~\ref{L-sl2-8} . This implies that $S$ has a least $n$ distinct zeros, a contradiction with Proposition~\ref{P-sl2-swf}. \hfill \qed \medskip

\end{add}%

\newpage
\bibliographystyle{plain}

\vspace{1cm}
\begin{flushleft}
PB: Universit\'{e} Grenoble Alpes and CNRS\\
Institut Fourier, CS 40700\\ 38058 Grenoble cedex 9, France\\
\verb+pierrehberard@gmail.com+\\[8pt]

BH: Laboratoire Jean Leray, Universit\'{e} de Nantes and CNRS\\
F44322 Nantes Cedex, France, and LMO, Universit\'e Paris-Sud\\
\verb+Bernard.Helffer@univ-nantes.fr+
\end{flushleft}%

\end{document}